\documentclass{article}
\usepackage{amsmath,amsthm,verbatim,amssymb,amsfonts,amscd, diagrams, graphics}
\theoremstyle{plain}

\theoremstyle{definition}

\theoremstyle{remark}

\newarrow{ul}---->
\newarrow{Backwards}<----

\newcommand{\into}{\hookrightarrow}
\newcommand{\Z}{\mathbb{Z}}

\newcommand{\R}{\mathbb{R}}
\newcommand{\C}{\mathbb{C}}
\newcommand{\bd}{\partial}

\def\co{\colon\thinspace}

\newarrow{Into}C--->

\begin{document}
\begin{center}
{\LARGE \textbf{Table of Contents for the}

\smallskip

\textbf{Handbook of Knot Theory}}
\end{center}
\smallskip

{\Large

\begin{center}
William W. Menasco and Morwen B. Thistlethwaite, Editors
\end{center}

\bigskip
\begin{enumerate}
\item Colin Adams, \emph{Hyperbolic knots}

\item Joan S. Birman and Tara Brendle, \emph{Braids: A survey}

\item John Etnyre, \emph{Legendrian and transversal knots}

\item Greg Friedman, \emph{Knot spinning}

\item Jim Hoste, \emph{The enumeration and classification of knots and links}

\item Louis Kauffman, \emph{Knot diagrammatics}

\item Charles Livingston, \emph{A survey of classical knot concordance}

\item Lee Rudolph, \emph{Knot theory of complex plane curves}

\item Marty Scharlemann, \emph{Thin position in the theory of classical knots}

\item Jeff Weeks, \emph{Computation of hyperbolic structures in knot theory}
\end{enumerate}
}

\pagebreak
\begin{center}
\textbf{\LARGE KNOT SPINNING}

\medskip

\textbf{\Large Greg Friedman}


\end{center}

\section{Introduction}

This exposition is intended to provide some introduction to higher-dimensional knots - embeddings of $S^{n-2}$ in $S^n$ - through spinning constructions. Once our shoe laces, those archetypal hand tools of knot theory, have been turned into spheres, how can we construct and visualize concrete examples of such knots? 

There are many important ways to construct higher-dimensional knots. If we are interested in \emph{algebraic knots}, we can look at the links of singularities of complex algebraic varieties in $\C^n$ (see, e.g., \cite{Miln, Du75}). We can also construct knots by \emph{surgery theory} (see \cite{LO} for a recent survey). There are powerful and complex tools for studying the knots that arise in these manners, but such construction methods  frequently do not allow one to ``see'' the knot. Often these knots can  be described only in terms of their algebraic invariants. We want to be able to visualize our knots, at least as far as it is possible to do so with our three-dimensional brains. This brings us to a series of constructions known as \emph{knot spinnings}. Many extensions have been made to Artin's original spinning technique, which dates back to 1925, but the various spinning constructions all have the appeal of being completely geometric in nature and thus highly visual. On top of providing a myriad of examples of importance in knot theory, these constructions provide an excellent introduction to thinking about higher-dimensional knots and higher-dimensional topology in general.

$\,$\linebreak

\medskip

Unfortunately, there do not seem to be many general references for knot theory in high dimensions, but we list a few sources that might be of interest for a beginner to the subject. More advanced references can be found in these sources. 

Colin Adams's popular treatment of knot theory in \emph{The Knot Book} \cite{Ad} contains  a chapter on visualizing high-dimensional knots. Dale Rolfsen's classic introduction to knot theory, \emph{Knots and Links} \cite{Rolf}, touches on the high-dimensional theory throughout, including spinning constructions on pages 85-87 and 96-99 and a ``Higher Dimensional Sampler'' in Chapter 11. The paper ``A survey of multidimensional knots'' by M. Kervaire and C. Weber in \cite{KW} provides a survey  through 1977. Andrew Ranicki's book \emph{High-dimensional Knot Theory: Algebraic Surgery in Codimension 2} \cite{Ran} provides a more modern and  extensive look at the theory from the point of view of algebraic surgery theory, while the article by Jerome Levine and Kent Orr \cite{LO} provides a more compact survey of high-dimensional knot theory via surgery. 

In addition, three recent books deal exclusively with knotted $2$-spheres (and other surfaces) in $\R^4$ and $S^4$. These are \emph{Braid and Knot Theory in Dimension Four} \cite{Kam} by Seiichi Kamada, \emph{Knotted Surface and Their Diagrams} \cite{CarSai} by J. Scott Carter and Masahico Saito, and \emph{Surfaces in $4$-Space} \cite{CKS} by Carter, Kamada, and Saito. All three of these books have a strong pictorial flavor, and each mentions knot spinning, the first book dealing with it more extensively in its Chapter 10 and the third book in Chapter 2.

\medskip

I thank Joan Doran for drawing the included figures. 

\section{Some basics}\label{S: basics}

\subsection{What is a knot?}
We  begin with the precise definition of a knot. 

Let $S^n$ be the \emph{$n$-dimensional sphere}, which we will be free to think of in several ways: as an abstract manifold, as the set of points in $\R^{n+1}$ unit distance from the origin, or as $\R^n$ compactified by adding a point at infinity. More generally, we will use $S^n$ to denote any object piecewise linearly (PL) homeomorphic to the sphere. Similarly, we will use $B^n$ to denote any object PL homeomorphic to the unit ball in $\R^n$, the set of points with distance $\leq 1$ from  the origin. The boundary of $B^n$, denoted $\bd B^n$, is PL homeomorphic to $S^{n-1}$.

With these conventions, a \emph{knot of dimension $n$} is a PL locally-flat embedding $K\co  S^{n-2}\into S^n$ or $K\co  S^{n-2}\into \R^n$. Recall that the piecewise linear (PL) condition simply means that there exist triangulations of $S^{n-2}$ and $S^n$ (or $\R^n$) with respect to which $K$ is simplicial, while $K$ is locally-flat if each point $K(x)$ in the image of $K$ has a neighborhood $U$  such that $(U, U\cap K(S^{n-2}))$ is PL homeomorphic to the standard coordinate pair $(\R^n, \R^{n-2})$. 
There is no real theoretical difference between letting $S^n$ or $\R^n$ serve as the codomain of the knot since we are free to rechoose  the point at infinity of $S^n$ so that the image of the knot  will lie in $\R^n\subset S^n$ (technically, we are replacing the knot with an \emph{equivalent} one; see below). Sticking with spheres has some technical advantages, and we will principally use spheres as the ambient space, though occasionally it will suit us to use $\R^n$ instead. 

\begin{figure}[!h]
\begin{center}
\scalebox{.6}{\includegraphics{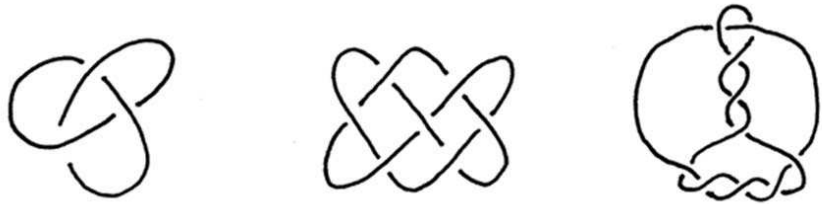}}
\end{center}
\caption{Some smooth  knots $S^1\subset \R^3$. A smooth knot can always be given a  PL structure by employing a suitable triangulation of $\R^3$ that contains the knot as a subcomplex.}
\end{figure}

\begin{figure}[!h]
\begin{center}
\scalebox{0.5}{\includegraphics{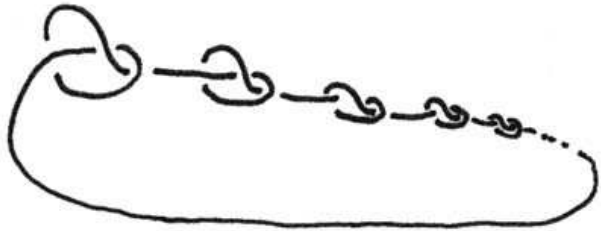}}
\end{center}
\caption{A ``wild'' knot that is not locally-flat at the point of infinite ``knottedness''.}\label{F: wild}
\end{figure}

The requirement that a knot be PL locally-flat is a common restriction, designed to avoid singularities of the embedding. For example, if we only required the embedding to be continuous,  ``infinite knottedness'' might occur (Figure \ref{F: wild}). Requiring the knot to be piecewise linear prevents this level of unpleasantness, but the local-flatness is also necessary to prevent other kinds of pathologies, such as local knotting,  that may occur when one is not working with differentiable maps - note that local-flatness certainly holds for differentiable embeddings by the Tubular Neighborhood Theorem. See \cite{GBF1} for a discussion of non-locally-flat knots.

In the classical dimension, $S^1\into S^3$, PL knots  and smooth knots are equivalent.
In fact, for any $n$,  any codimension two PL locally-flat embedding can be made differentiable (see \cite[Corollary 6.8]{RS68}).  However, in high dimensions the smooth structure on the embedded sphere may not be the standard one. By working in the PL category, we allow these knots but do not concern ourselves with any eccentricities in their smooth structures. All this  being said, most of the constructions we will discuss work equally well in the differentiable category and generate knots with the standard smooth structure provided we start with a knot with the standard smooth structure.  
Knot spinning can also be done on non-locally-flat PL knots provided some minor extra care is employed  (see \cite[\S 4.3]{GBF1}) or with topological non-locally-flat knots provided that the embedding is flat at some point. In the sequel, we will stick with PL locally-flat knots for convenience and consistency.

\subsection{Knot equivalence}

Now, if you have a knotted string lying on your desk and you pick it up and move it someplace else, we would like to think of it as the same knot. Thus we should really consider equivalence classes of knots. We call two knots $K_0, K_1\co  S^{n-2}\to S^n$ \emph{equivalent} if there is an orientation-preserving PL homeomorphism $f\co S^n\to S^n$ such that $fK_0(S^{n-2})=K_1(S^{n-2})$. In other words, $f$ should take the image of  $K_0$ to the image of $K_1$. In particular, this will be true if there is a PL ambient isotopy of $S^n$ taking $K_0(S^{n-2})$ to $K_1(S^{n-2})$. In fact, this stronger condition is sometimes used as the definition of knot equivalence. Since we work in the PL-locally-flat category, these two conditions are equivalent (see \cite[Proposition 1.10]{BZ}), but the analogous equivalence does not hold in the smooth category due to the failure of the Alexander trick.

It is a standard abuse, in which we shall engage freely, to use the word ``knot'' and the same symbol, $K$, to refer to the equivalence class of the knot $K$ or even to the image of $K$. 

We refer to  $K\co S^{n-2}\into S^n$ as \emph{$n$-dimensional} or an \emph{$n$-knot}. This is    not a universal notation; it is perhaps more standard to refer to such a knot as an $n-2$ knot. We also refer to the knots $K\co  S^1\into S^3$ as \emph{classical knots}.

One also sometimes speaks of \emph{oriented equivalence} for which $S^{n-2}$ is given a fixed orientation and it is required that the orientation-preserving PL homeomorphism $f\co S^n\to S^n$
taking $K_0(S^{n-2})$ to $K_1(S^{n-2})$ also preserves the orientation of these subspaces. However, we will not impose this stricter condition except when stated explicitly.

\subsection{The unknot and toroidal decompositions of $S^n$}

The \emph{unknot} in dimension $n$ is the equivalence class of the ``standard embedding'' $S^{n-2}\subset S^n$. In other words, if $S^n=\{\vec x\in \R^{n+1}\mid |\vec x|=1\}$, then the unknot can be represented as $\{\vec x\in \R^{n+1}\mid |\vec x|=1, x_{n+1}=x_n=0\}$. The classical $3$-dimensional unknot is equivalent to the unit circle in the $x$-$y$ plane in (the compactified) $\R^3$. 

By setting more coordinates equal to $0$, we can define standard embeddings of any sphere into any other sphere of higher dimension. This leads to  nice decompositions of $S^n$ into two generalized solid tori: The standard embedding of $S^m$ into $S^n$, $m<n$, has a tubular neighborhood  PL homeomorphic to $S^m\times B^{n-m}$, and the complement of the interior of this neighborhood is PL homeomorphic to $B^{m+1}\times S^{n-m-1}$. Thus $S^n$ can be decomposed as a union of  $S^m\times B^{n-m}$ and $B^{m+1}\times S^{n-m-1}$,  identified in the obvious way along their common boundary $S^m\times S^{n-m-1}$. $$S^n=S^m\times B^{n-m}\bigcup_{S^m\times S^{n-m-1}} B^{m+1}\times S^{n-m-1} $$ This follows, e.g., from the fact that $S^n$ can be written as the join of $S^m*S^{n-m-1}$. 
The most familiar case is the standard genus one Heegard decomposition of $S^3$, in which a neighborhood of the unknot $S^1\subset S^3$ is the  solid torus $S^1\times B^2$, whose complement is another solid torus, $B^2\times S^1$.

\begin{figure}[!h]
\begin{center}
\scalebox{0.6}{\includegraphics{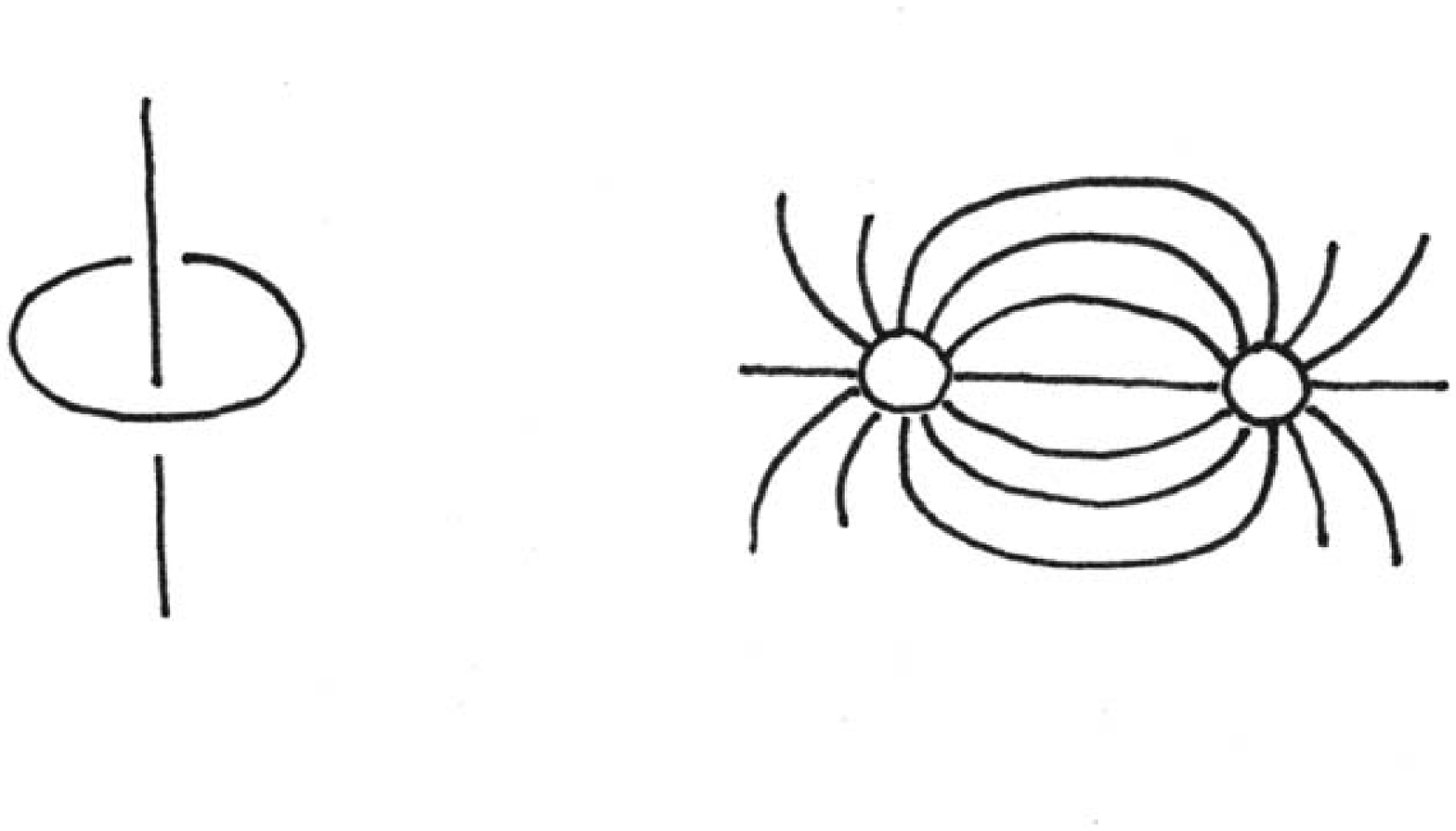}}
\end{center}
\caption{The decomposition of $S^3$ (thought of as $\R^3$ plus a ``point at infinity'') into two solid tori $S^1\times B^2$ and $B^2\times S^1$. The lefthand picture shows the circular cores of the tori (the vertical line becomes a circle as it wraps through the point at infinity). The righthand picture shows a slice along the  $y$-$z$ plane: The two disks are a slice of one solid torus (cut a donut in half and then view it on end), while the arcs represent slices of the meridional disks  of the other solid torus.}\label{F: torus decompostion}
\end{figure}

\subsection{A useful excision}\label{S: ex}

\begin{figure}[!h]
\begin{center}
\scalebox{0.5}{\includegraphics{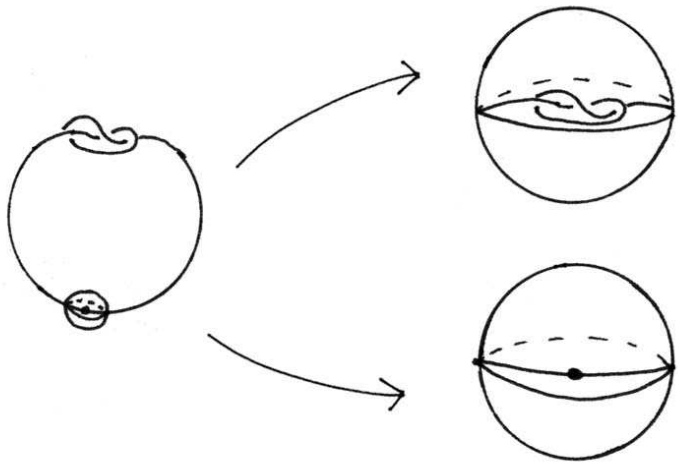}}
\end{center}
\caption{Removing a trivial neighborhood from a knotted circle to obtain a knotted arc.}\label{F: complement}
\end{figure}

We conclude this introductory section with one other construction that will be used repeatedly. Consider an $n$-dimensional knot $K$, and choose any  point $x\in K\subset S^n$ (here - by our standard notational abuse - we use ``$K$'' to represent the image of the knot).  Since $K$ is PL locally-flat, there is a neighborhood $B^n_{-}$ of $x$ in $S^n$ such that $(B^n_{-}, B^{n-2}_{-}):=(B^n_{-}, B^n_{-}\cap K)$ is  PL homeomorphic to an unknotted ball pair, i.e. it is PL homeomorphic to the standard ball pair $(B^n, B^n\cap \R^{n-2})$ in $\R^n$.  
Since the closure of the complement of a PL $n$-ball in $S^n$  is also an $n$-ball, the closures of the complements  $S^n-B_{-}^n$ and $K-B_{-}^{n-2}$ will each be balls, and we label this complementary pair by $(B_K^n, B_K^{n-2})$. The ball $B_K^{n-2}$ may be knotted in $B_K^n$ (Figure \ref{F: complement}).

We also observe that the common boundary of the pairs $(B^n_{-}, B^{n-2}_{-})$ and $(B_K^n, B_K^{n-2})$ is the \emph{unknotted} pair of spheres $(S^{n-1}, S^{n-3})$ (since it is PL homeomorphic to the boundary of the standard ball pair). In what follows, it will often be convenient to identify this with the \emph{standard unknot}, which we have already discussed.

What if we choose the neighborhood of a different point to remove in this construction? It turns out that we get the same pair $(B_K^n, B_K^{n-2})$ up to PL homeomorphism. To see this, consider the ball neighborhoods of two different points. We can  simply slide one ball to the other along the knot, which complementarily takes the complement of one neighborhood to the complement of the other. Note that while this idea has nice intuitive appeal, it does require some technical checking to ensure that such sliding is always allowed. However, this theory is well-established, and we avoid going too far afield to visit the details here (see, e.g., Chapter 6 of Hudson \cite{Hud}).

\section{Basic spinnings}\label{S: basic}

\subsection{Simple spinning}\label{SS: spin}

As we know, there are an infinite number of  knots $S^1\into S^3$  and myriad examples can be created by anyone with a piece of string and some time on their hands (classifying these knots is another matter!). To get knots of higher dimensions requires a little bit more ingenuity. One method is to get high-dimensional knots from knots of lower dimension by spinning them. 

The earliest spinning construction is due to Emil Artin in 1925 \cite{Ar}. Artin used spinning to construct $4$-dimensional knots from classical knots, but the same idea can be used to create an $n+1$ dimensional knot from any $n$-dimensional knot. This construction is generally referred to just as ``spinning'', but we will call it \emph{simple spinning} to differentiate it from the more general constructions to follow. 

In this section, it will be most convenient to consider knots in $\R^n$ instead of $S^n$ (see Section \ref{S: basics}), though of course we can easily transform from one type to the other by adding or removing a point at infinity. 

\begin{figure}[!h]
\begin{center}
\scalebox{.5}{\includegraphics{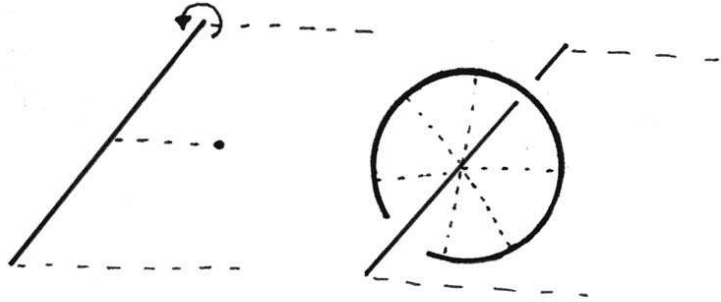}}
\end{center}
\caption{Spinning a point in the half-plane around the axis.}\label{F: rotate}
\end{figure}

 To see the basic idea, consider  the upper half plane $H^2=\{(x,y)\in \R^2\mid y\geq 0\}$ and choose a point $(x_0,y_0)\in H^2$ with $y_0>0$. Now rotate $H^2$ around the $x$-axis in $\R^3$. The point will sweep out a circle (Figure \ref{F: rotate}). Analytically, the circle will be parametrized in $\R^3$ by the set of points $(x_0, y_0\cos\theta, y_0 \sin\theta)$, as $\theta$ runs from $0$ to $2\pi$ (assuming that we rotate counterclockwise as seen from the positive $x$-axis looking in the negative $x$ direction). 

\begin{figure}[!h]
\begin{center}
\scalebox{.5}{\includegraphics{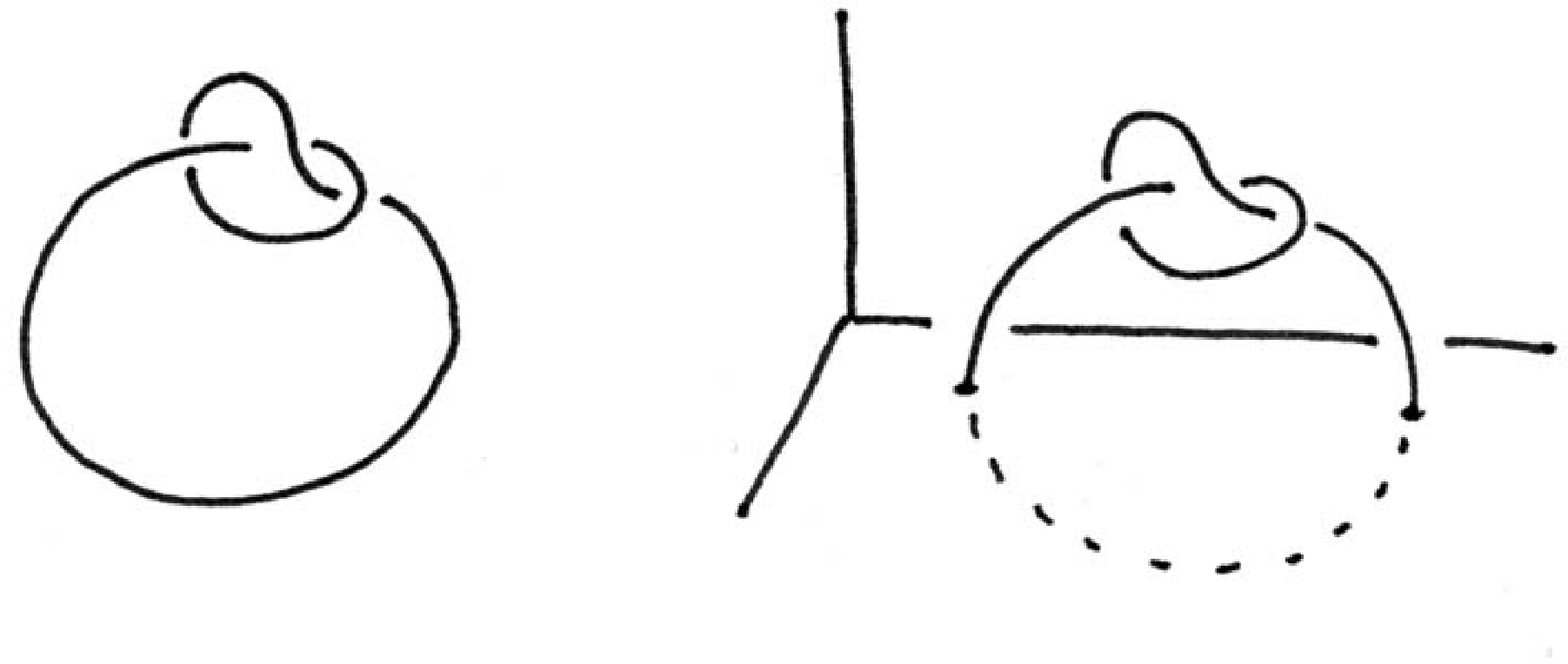}}
\end{center}
\caption{Turning a knotted circle into a knotted arc in the upper half space in order to spin it about the plane.}\label{F: simple spin}
\end{figure}

\begin{figure}[!h]
\begin{center}
\scalebox{.5}{\includegraphics{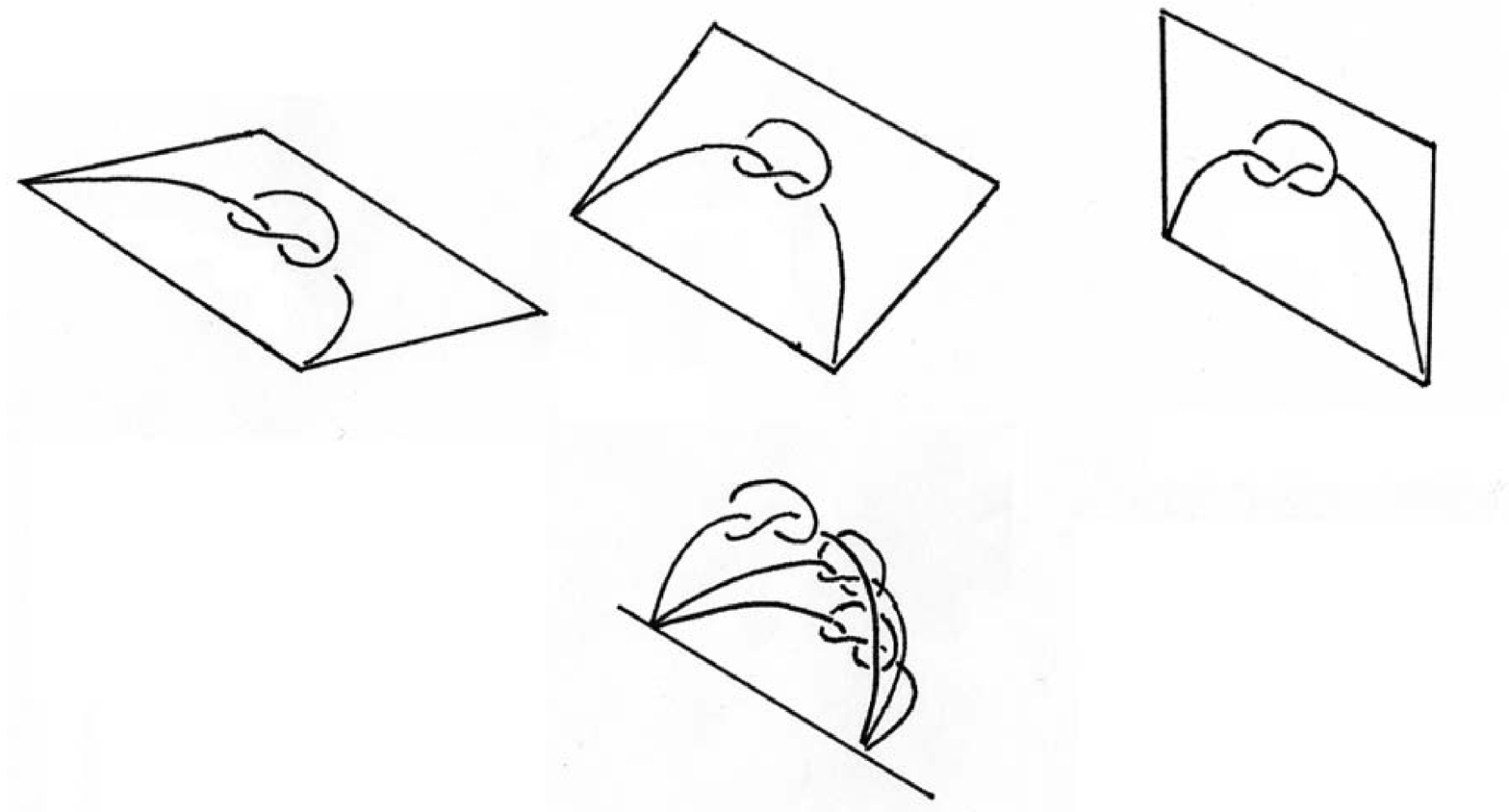}}
\end{center}
\caption{A schematic of knot spinning }\label{F: knotfilm}
\end{figure}

\begin{figure}[!h]
\begin{center}
\scalebox{.5}{\includegraphics{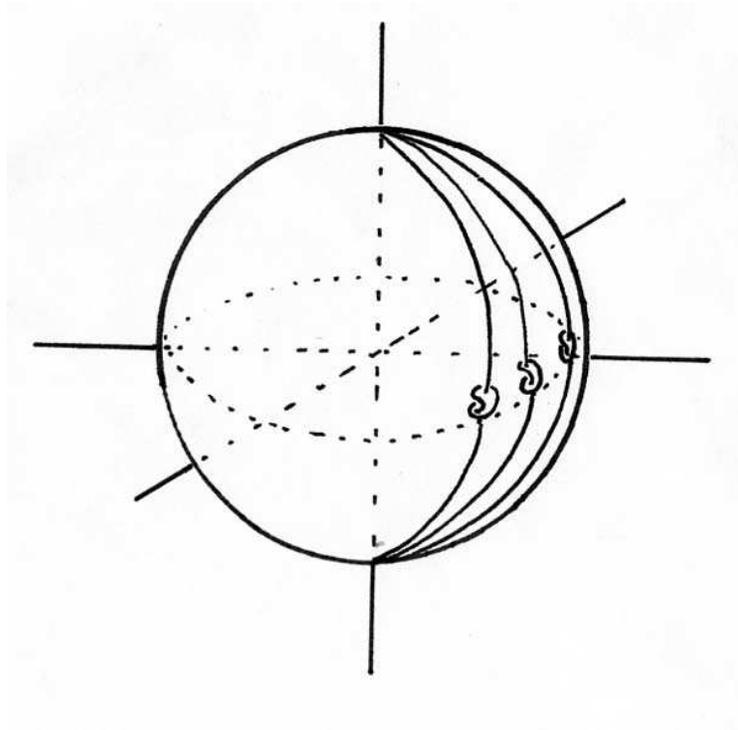}}
\end{center}
\caption{A second schematic of knot spinning. }\label{F: rotating}
\end{figure}

To see how this applies to knots, let us consider a knot $K$ in $\R^3$. Up to equivalence, we can arrange for $K$ to lie in the upper half space $H^3=\{(x,y,z)\mid z\geq 0\}$ except for an unknotted arc that dips below the $x$-$y$ plane $\R^2=\{(x,y,z)\mid z=0\}$ (Figure \ref{F: simple spin}). Let us remove the interior of this unknotted arc; what remains is a knotted arc in $H^3$ with its endpoints (and only its endpoints) in $\R^2$. We can now rotate $H^3$ around $\R^2$ in $\R^4$ just as we rotated $H^2$ around $\R^1$ in $\R^3$. Analytically, we parametrize by $\theta$, and each point $(x,y,z)$ in the upper half space sweeps out the circle $(x,y,z\cos\theta, z\sin \theta)$. Note that $\R^2$ remains fixed. By thinking about how the longitude lines swing around the globe with the north and south poles remaining fixed, we can imagine how the knotted arc gets spun into the image of a 2-sphere $S^2$. Thus, by spinning, we obtain a knotted $S^2$ in $\R^4$ (Figures \ref{F: knotfilm} and \ref{F: rotating}).

You might be asking, what if we had chosen a different way to split
our original knot into a knotted arc? It turns out that we get the same spun
knot, essentially for the same reason by which we noted in Section \ref{S: ex} that $B_K^n$ is independent of the choice of $B^n_-$.  In fact, notice that if we start with our knot in $S^3$, then
our construction to get a knotted arc in the upper half space by removing an unknotted arc in the lower half space is
completely equivalent to the construction of $B_K^n$ by removing a small ball neighborhood of a point on the knot.

This simple spinning construction already has several important ramifications. For example, it can be shown very easily  (see, e.g., \cite{Rolf}) that the fundamental group of the complement of this spun knot in $\R^4$ (its knot group) is isomorphic to the knot group of  our original knot in $\R^3$. Based on known results about knots in $\R^3$, this  implies the existence of an infinite number of inequivalent knots in $\R^4$. By contrast, Dennis Roseman  showed in \cite{Ro75} that the spins of two distinct knots may be equivalent. For example, he showed that spinning the square knot yields the same $4$-knot as does spinning the granny knot. We discuss a much stronger result along these lines at the end of Section \ref{S: superspin}.

The construction for higher dimensions is similar. We begin with a knot $K\co  S^{n-2}\to \R^n$. Again, we can manipulate the knot within its equivalence class so that it lies mostly in the upper half space $H^n=\{(x_1,\ldots, x_n)\in \R^n\mid x_n\geq 0\}$ and so that the intersection of the knot with the lower half space is  unknotted. We then remove the interior of this unknotted ball to obtain the complementary  knotted ball $B^{n-2}$ in $H^n$. Its intersection with $\R^{n-1}$ is an  unknotted $S^{n-3}$. Now we spin $H^n$ around $\R^{n-1}$ into $\R^{n+1}$ so that each point $(x_1,\cdots, x_n)$ sweeps out the circle $(x_1,\cdots, x_{n-1},x_n \cos \theta, x_n\sin \theta)$. 

It is a little harder now to see that our knotted ball in the upper half plane gets spun into a sphere $S^{n-1}$, but the idea of pivoting a longitude around its poles extends to higher dimensions. To see this, we employ the following coordinate analysis, which will also be useful when we need to describe more general spinnings, below.

Consider $S^{n-1}$ as the unit sphere in $\R^n$,  $S^{n-1}=\{\vec x\in \R^n\mid |\vec x|=1\}$, and consider $\R^n$ as $\R^{n-2}\times \R^2$. Then we can
define the latitude for a point $y\in S^{n-1}$ as its orthogonal projection onto the
$\R^{n-2}$ factor and the longitude of $y$ as the angular polar coordinate of its
projection onto the $\R^2$ factor. Hence the latitude is always
well-defined, while the longitude is either undefined or a unique angle, dependent upon whether or not $y$ lies in $\R^{n-2}\times 0$. Notice that in the
case where the longitude is undefined, the point on the sphere is uniquely
determined by its latitude (just as on a globe). To simplify the notation in abstract cases, we will simply refer to the
latitude-longitude coordinates $(z, \theta)$, whether $\theta$ is defined or not. Then the point $(z,\theta)$ in $S^{n-1}$ corresponds to the point in $\R^n$ determined by the rectangular coordinates $(z, r\cos\theta, r\sin\theta)$ for $z\in \R^{n-2}$, $r\geq0$, $\theta\in [0,2\pi)$, and such that $|z|^2+r^2=1$ (note that this determines $r$, given $z$). 

Now, consider the closure of the set of points in $S^{n-1}$ with  fixed longitude $\theta=0$. These points can be written in rectangular coordinates as $(z, r, 0)$, with $|z|^2+r^2=1$ and $r\geq 0$. This set is homeomorphic to a ball $B^{n-2}$ (in fact, it is the graph of $r=\sqrt{1-|z|^2}$). Its boundary is the $n-3$ sphere with $|z|^2=1$ and $r=0$, and we call this  boundary the \emph{generalized pole} of $S^{n-1}$. Now for each point $(z,r,0)$ in rectangular coordinates, we can spin to get the set of points $(z,r\cos \theta, r\sin \theta)$ as $\theta$ runs from $0$ to $2\pi$. The points $(z,0,0)$ of the generalized pole remain fixed, and the rest of the $0$ longitude sweeps out the rest of the sphere. 
Analogously, as we spin a knot, the knotted ball  sweeps out a knotted sphere. We leave it to the reader to formulate a precise analytic description.

Preservation of  knot groups under simple spinning continues to hold in this higher-dimensional setting, and by iterating the spin construction, we establish the existence of an infinite number of inequivalent knots in any dimension $n\geq 3$.

\subsection{Superspinning}\label{S: superspin}

Having spun knots in  circles, how about spinning around higher dimensional spheres? \emph{Superspinning} of classical knots was introduced by E.C. Zeeman \cite{Z60} and D.B.A. Epstein \cite{Ep} separately in 1960. Zeeman used the spinning of classical knots about  spheres to show that it is possible to embed two $n-2$ spheres in $S^n$ such that each sphere in unknotted but the pair is linked, i.e. there is no homeomorphism of $S^n$ sending one knot to the ``northern hemisphere'' and the other knot to the ``southern hemisphere''. Epstein strengthened these results to show that  
two $n-2$ spheres can be embedded in euclidean $n$-space  in each of the following ways: (i) neither can be shrunk to a point in the complement of the other; (ii) one can and one cannot be shrunk to a point in the complement of the other.
In  1970, Sylvain Cappell \cite{C70} generalized this construction as a way to construct an $n+p$ dimensional knot from any $n$-dimensional knot by spinning it around a $p$-sphere $S^p$, $p\geq 1$. He called this method ``superspinning'', though it is sometimes referred to in the literature as $p$-spinning. Cappell utilized superspinning to demonstrate the existence of knots whose complements are homotopy equivalent but not homeomorphic. 

This time let us jump straight to the general construction, $p$-spinning for a knot $K\co S^{n-2}\into S^n$. 

Imagine $ S^p \times B^n$ embedded in  our standard unknotted way in $S^{p+n}$ so that  we can write $$S^{p+n}=[ S^p \times B^n]\bigcup_{S^p\times S^{n-1}} [B^{p+1}\times S^{n-1}]$$ (see Section \ref{S: basics}). Here $\bigcup_{S^p\times S^{n-1}} $ indicates that we are gluing the two spaces along their common boundary  $S^{p}\times S^{n-1}$.   We can decompose the unknot $S^{p+n-2}$ in $S^{p+n}$ by its intersections with the pieces of  this decomposition  as $$[ S^p \times B^{n-2}]\bigcup_{S^p\times S^{n-3}} [B^{p+1}\times S^{n-3}].$$ 
\begin{diagram}
S^p\times B^{n-2}&\cup & B^{p+1}\times S^{n-3}\\
\dInto&&\dInto\\
S^p\times B^{n}&\cup & B^{p+1}\times S^{n-1}
\end{diagram}
(Here we think of $B^{n-2}$ as the unknotted subset of $B^n$ given by setting the last two coordinates to $0$.)
This is one of our standard decompositions of $S^{p+n-2}$, but now we see it lying within a decomposition of the larger sphere $S^{p+n}$. We can write the pair of spaces more compactly as $$S^p\times (B^n, B^{n-2})\cup B^{p+1}\times (S^{n-1},S^{n-3}),$$ and we should think of this as the product of $S^p$ with a trivial (unknotted) ball pair, ``capped off'' by another standard piece.  

Try to picture this decomposition of the unknotted $S^2$ in $S^4$, taking $p=1$, $n=3$ and recalling that $S^0$ is a pair of points. In this case, $S^2$ decomposes into a neighborhood of the equator and neighborhoods of the north and south poles. The decomposition of $S^4$ consists of a neighborhood in $S^4$ of  the equator of $S^2$  and its complement, which is $B^2\times S^2$. So, written as pairs, we have the unknot $(S^{4},S^{2})$ decomposed into  $S^1\times (B^3,B^1)$, which will play the important role in our spinning construction, and $B^2\times (S^2,S^0)$, ``the rest''. 

Now, within this construction, we can replace  each trivial pair $(B^n, B^{n-2})$ in the product $S^p\times (B^n, B^{n-2})$  with the knotted ball pair  $(B_K^n, B_K^{n-2})$, obtained from $K$ as in Section \ref{S: ex}. In other words, we construct $$S^p\times (B^n_K,B^{n-2}_K)\bigcup B^{p+1}\times (S^{n-1},S^{n-3}) ,$$ and we define the superspun knot $K^*$ to be the subset given by  $$[S^p\times  B_K^{n-2}]\bigcup_{S^p\times S^{n-3}} [B^{p+1}\times S^{n-3}].$$ As $K^*$ is PL homeomorphic to the standard decomposition of $S^{n+p-2}$, we see that $K^*$ is  a sphere of dimension $n+p-2$ knotted in $S^{n+p}$. 

If $p=1$, superspinning $K$ gives us the same simple spun knot that we obtained in the Section \ref{SS: spin}. Can you see why? Try thinking about $1$-spinning classical knots. 

It turns out that the knot group of a superspun knot is also the same as the knot group of the original knot, but in general, superspinning does not create the same knots as does iterated simple spinning. 

Also as for simple spinning, $p>1$ superspinning may take inequivalent knots to equivalent knots: Cameron Gordon showed in \cite{Go76} that all superspun knots are amphicheiral, i.e., they are oriented equivalent to the knot obtained by reversing the orientation of both $S^{n+p}$ and $K^*$ (this is sometimes called $(-)$-amphicheirality). As a corollary, this generalizes Roseman's result, cited above, on the equivalence of the spun granny knot with the spun square knot, and it implies that the $p>1$ superspins of inequivalent knots may be equivalent. Another result along these lines was obtained by Cherry Kearton, who showed in \cite{Kea91} that the superspins of two classical knots $K_1, K_2\subset   S^3$ are equivalent if and only if their knot groups are isomorphic. This is false, however, for the superspins of knots of higher dimension (see, e.g., \cite{Kea84}).

\subsection{Frame spinning}

Even more general than superspinning is \emph{frame spinning}: why limit ourselves to spinning about spheres? How about other manifolds? Frame spinning was introduced by Dennis Roseman in 1989 \cite{Ro89}, though the name is due to Alexander Suciu \cite{Su92}, who used frame spinning to construct new examples of inequivalent knots that have the same complement (this is a phenomenon that occurs only for knots above the classical dimension, though at most two higher-dimensional knots can share a given complement; see, e.g., \cite{LS}).

To describe frame spinning, let us once again begin with an $n$-dimensional knot $K$. This time, however, our additional data comes in the form of an $m$-dimensional manifold $M^m$ embedded in $S^{n+m-2}$ with a framing $\phi$. This last condition means that we in fact consider an embedding $\phi\co  M^m\times B^{n-2}\into S^{n+m-2}$. Furthermore, we assume that $S^{n+m-2}$ is embedded in the standard, unknotted way into $S^{n+m}$ with the standard framing as in the generalized torus decomposition.  Putting these framings together, we get a pair of tubular neighborhoods of $M^m$ in  $(S^{n+m},S^{n+m-2}) $ of the form $N=M^m\times (B^n, B^{n-2})$, where each $(B^n, B^{n-2})$ is an \emph{unknotted} ball pair (although the exact embedding of $N$ into $S^{n+m}$ depends on our choice of framing $\phi$). 

\begin{figure}[!h]
\begin{center}
\scalebox{.6}{\includegraphics{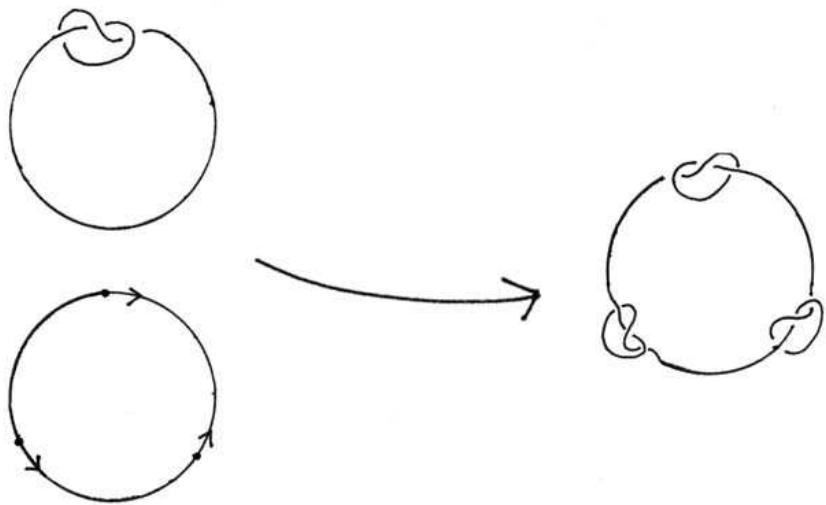}}
\end{center}
\caption{The trefoil knot spun about the manifold $M$ consisting of three disjoint points in $S^1$. Note that the framing at each point (indicated by an arrow that depicts the orientation of the framing) determines how to attach the knot.}
\end{figure}

The idea now is to take all of those unknotted ball pair fibers in $N$ and replace them with our knotted ball pair $(B_K^n, B_K^{n-2})$ as we did for superspinning. In other words, having used the framing to identify the neighborhood pair $N$ as $M^m\times (B^n, B^{n-2})$, we remove it, and then replace it with $M^m\times (B_K^n, B_K^{n-2})$, glued in along the same framing. Thus, our frame spun knot will be $$(S^{n+m-2}-M^m\times B^{n-2})\bigcup_{M^m\times S^{n-3}}M^m\times B_K^{n-2}$$ sitting inside  the $n+m$ sphere $$(S^{n+m}-M^m\times B^{n})\bigcup_{M^m\times S^{n-1}}M^m\times B_K^{n}.$$ In the special case where $M^m$ is the sphere $S^m$ embedded in the standard way and with standard framing in $S^{n+m-2}$, we recover $m$-superspinning (why?).

If the manifold $M$ has multiple components, or even components of different dimensions, then we can spin different knots (also possibly of different dimensions) around each component. It is possible to generalize this construction even further, but first we should study some other types of spinning.

\section{Spinning with a twist}
\subsection{Twist spinning}
\emph{Twist spinning}, introduced by E.C. Zeeman in 1965 \cite{Z65}, was an early generalization of Artin's simple spinning construction. Again, we begin with an $n$-dimensional knot and obtain an $n+1$ dimensional knot, but the difference between simple spinning and twist spinning can be illustrated celestially: As the moon orbits the Earth, it always keeps the same face towards the Earth. This is analogous to simple spinning in which the knot is rotated around the plane but always keeps ``the same face'' towards the plane serving as the axis of rotation. By comparison, twist spinning is like the Earth orbiting the sun: as the earth orbits, it also rotates around its own axis. 

\begin{figure}[!h]
\begin{center}
\scalebox{0.5}{\includegraphics{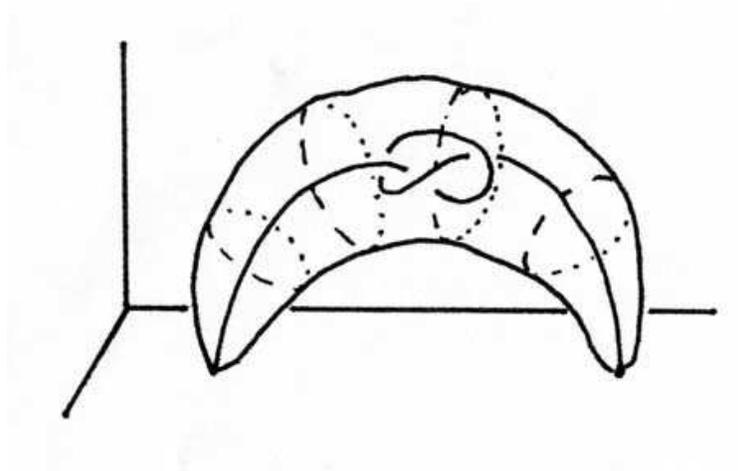}}
\end{center}
\caption{We rotate a ball around the knotted arc.}\label{F: banana}
\end{figure}

\begin{figure}[!h]
\begin{center}
\scalebox{0.5}{\includegraphics{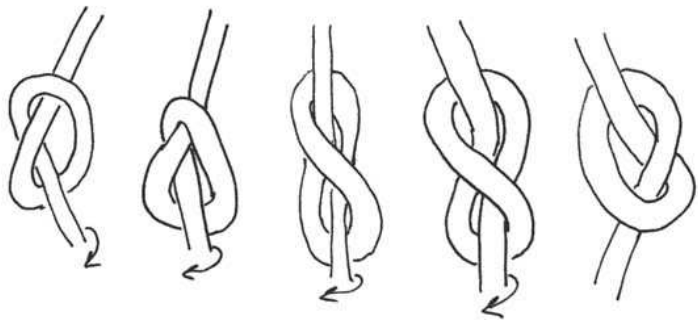}}
\end{center}
\caption{A 180 degree twist of the  trefoil knot (thickened for improved visualization).}\label{F: spin}
\end{figure}

Before giving a general formula, let us consider heuristically the case of twist spinning a classical knot. As in the simple spinning construction, we replace the knot with a knotted arc in the upper half space whose endpoints lie in the $x$-$y$ plane. We can also assume that the knotted part of the arc is contained within a ball whose intersection with the arc is its north and south poles (Figure \ref{F: banana}).  Now, as we rotate half-space around the plane as in simple spinning, we simultaneous spin this ball on its axis (Figure \ref{F: spin}). It is only necessary that the end result lines up with the starting position, so we are free to spin the ball on its axis any integral number $k$ times as we rotate $H^3$. 

Let us be more specific. Given an $n$-knot $K$, then just as for superspinning about $S^1$ (which is equivalent to simple spinning), we decompose the $n+1$ dimensional unknot as  the two space pairs $S^1 \times (B^n, B^{n-2})$ and $B^2\times (S^{n-1}, S^{n-3})$. To superspin, we simply removed $S^1\times (B^n, B^{n-2})$ and glued in $S^1\times (B_K^n, B_K^{n-2})$, reattaching along the original boundary $S^1\times(S^{n-1},S^{n-3})$.  In order to create the $k$-twist spin, however, we glue  in the following way: we represent points in the common boundary $S^1\times S^{n-1}$ by  $(\eta, z, \theta)$, where $\eta\in S^1$ and $(z, \theta)$ are latitude-longitude coordinates for $S^{n-1}$ such that the unknotted $S^{n-3}\subset S^{n-1}$ is the generalized pole (see Section \ref{SS: spin}).  If $(\eta, z,\theta)$ is such a point in the boundary of $B^2\times S^{n-1}$, we attach that point to the boundary of $S^1 \times B_K^n$ by $(\eta, z,\theta)\to (\eta, z, \theta+k\eta)$. The addition here is standard angle addition in the circle, which we can think of as $\R/2\pi\Z$.  
In this way, as we glue the pieces together, we introduce a $k$-fold twisting by rotations of the longitude coordinate. You should convince yourself that this procedure corresponds to our earlier heuristic description. 

\begin{figure}[!h]
\begin{center}
\scalebox{.5}{\includegraphics{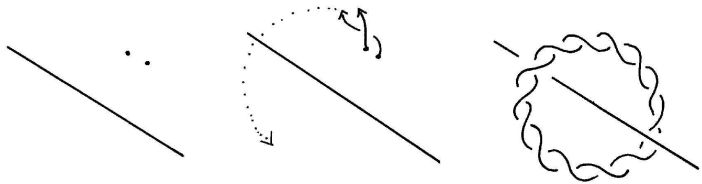}}
\end{center}
\caption{A low-dimensional schematic of twist-spinning.}\label{F: twistschem}
\end{figure}

If your first instinct is to think that the twisting doesn't add anything to the spinning since  the knotted arc lands back where it started, you should consider the following toy example. Recall from Section  \ref{SS: spin}
our original toy example of simple spinning in which we rotated a point in the upper half plane around the $x$-axis. This time, however, imagine two points in the upper half plane. As we sweep the half plane around the $x$-axis, let these two points rotate around each other in the plane $k$ times, where $k$ is any integer. At the end of the process, we have two curves in space that link each other $k$ times (Figure \ref{F: twistschem}. This isn't quite the same procedure as twist spinning, but it should illustrate the idea that interesting things can happen if we deform as we spin.

Zeeman showed in \cite{Z65} that a twist spun knot depends only on $K$ and $|k|$, i.e. $k$-twist spinning and $-k$-twist spinning yield the same knot. Furthermore, he proved the slightly surprising fact that any $1$-twist spun knot (and hence also any $-1$-twist spun knot) is  unknotted! This is actually a corollary of the much stronger theorem in the same paper stating that any $k$-twist spun knot is a fibered knot with fiber the punctured $k$-fold branched cover of $S^n$ determined by the knot being spun. It has also been shown, by Cameron Gordon in his thesis \cite{Go}, that if $k$ and $l$ are coprime, then a $k$-twist spin followed by an $l$-twist spin yields the unknot; this was generalized by Tristram (unpublished), who showed that any sequence of $k_i$ twist spins, where the $k_i$ are coprime, is unknotted. A short proof of both statements can be found in \cite{Kea72}.  Deborah Goldsmith and Louis Kauffman found another generalization by showing in  \cite{GK78} that, if $L_{k,l}(K)$ indicates the $l$-twist spin of the $k$-twist spin of $K$, then  $L_{k,l}(K)$ is equivalent to one of $L_{0,g}(K)$ or $L_{g,g}(K)$, where $g=\text{g.c.d.}(k,l)$.

Unlike the constructions of Section \ref{S: basic}, twist spinning does not preserve knot groups. In fact, if $G$ is the knot group of  $K$, then the  group of the  $k$-twist spin of $K$ is isomorphic to $(\Z\times G)/\langle t^{-1}g^k \rangle$, where $t$ is a generator of $\Z$ and $g\in G$ represents a meridian of the knot $K$ (see, e.g., 
\cite{GBF4}). Note that if $k=0$, i.e. we simple spin, this group is just $G$. It is less easy to see, though it must be true by Zeeman's theorem, that if  $k=\pm 1$, then this group is $\Z$, the knot group of the unknot! Actually, this is not hard to show algebraically, using the fact that any knot group is the normal closure of any element representing a meridian.

\subsection{Frame twist spinning}

Now that we have seen how to add twisting to Artin's basic spinning construction, can we add some kind of twisting to our other spins? For superspinning about spheres of dimension greater than $1$, the answer is no! This is because $\pi_n(S^1)=0$ for all $n>1$, which implies that we cannot use higher dimensional spheres to parametrize spinning. Any attempt at twisting can be deformed to give back ordinary superspinning. On the other hand, since $\pi_1(S^1)=\Z$, there are countably many maps $S^1\to S^1$ that cannot be so deformed, and the element $k$ in $\pi_1(S^1)$ corresponds to $k$-twist spinning. 

However, where superspinning fails to be twistable, frame spinning does allow a twist if the manifold $M^m$  admits a map $M^m\to S^1$ that cannot be deformed into the trivial map to a point. Just as for twist spinning, this map provides us with enough data to alter the gluing map of the construction by twisting the longitude coordinate of $B_K^n$ as we glue. The gluing of the latitude coordinates is once again controlled by the framing $\phi$ of $M$. 

So let us be specific. Recall that, in frame spinning, we used the framing of $M^m$ in $S^{n+m-2}$ along with the trivial framing of $S^{n+m-2}$ in $S^{n+m}$  to identify a neighborhood $N$ of $M^m$ in $(S^{n+m}, S^{n+m-2})$ with the product $M^m\times (B^n, B^{n-2})$. Then we replaced $M\times (B^n, B^{n-2})$ with  $M\times (B_K^n, B_K^{n-2})$ and glued it back in \emph{along the same framing}. Suppose, however, that we are given a map $\tau\co  M^m\to S^1$. Then we can use this map to augment the gluing with a twist along the longitude. This is done as follows: we use the framings to assign coordinates $(x,z,\theta)$ to the boundary $M\times (S^{n-1},S^{n-3})$ of the neighborhood $N$.  Here $x\in M$ and $(z, \theta)$ are latitude/longitude coordinates on $S^{n-1}$ such that the unknotted $S^{n-3}$ is the generalized pole. The boundary of the complement of $N$ in $S^{n+m}$ possesses the same coordinates, as   these two boundaries agree.
Again, we cut out $N$ and replace it with $M^m\times (B_K^n, B_K^{n-2})$, which we glue back in, but instead of gluing the point $(x,z,\theta)$ in $\bd (S^{n+m}-N)$ right back to its counterpart in $\bd N$, we glue it by the attaching map $f\co  (x,z,\theta)\to (x,z,\theta+\tau(x))$, where again the addition in $S^1$ is standard angle addition. 

In other words, we form 
$$
[(S^{n+m}, S^{n+m-2})-M^m\times (B^{n},B^{n-2})]
\bigcup_f [M^m\times (B_K^{n}, B_K^{n-2})],
$$
where $\bigcup_f$ indicates a gluing via the attaching map given above.

This construction was introduced by the author in his dissertation. He goes on to calculate various algebraic invariants of frame twist-spun knots, including Alexander polynomials in \cite{GBF1} and the knot groups in \cite{GBF4}. If the map $\tau$ is homotopic to the trivial map, we recover Roseman's frame spinning. If $M=S^1$ and $\tau$ is the map that wraps the circle around itself $k$ times, we recover Zeeman's $k$-twist spinning. 

\section{More general spinnings}

\subsection{Deform spinning}

An even more general class of spinning constructions is known to exist. The first example, \emph{roll spinning}, was introduced in a short paper by Ralph Fox in 1966 \cite{Fox}. Unfortunately, Fox provided only one example and did not include specific details of the construction, which has led to some contention over the exact definition of roll spinning. In 1979, R.A. Litherland \cite{Lith}  provided a formal definition and a generalization, \emph{deform spinning}, of which  both roll spinning and  twist spinning are special cases. According to Masakazu Teragaito \cite{Ter}, Fox's original construction is actually an example of what Litherland calls \emph{symmetry spinning}, but by now Litherland's definition of roll-spinning is the one that has caught on.

Deform spinning is another construction that takes $n$-knots to $n+1$-knots.

The tersest description of deform spinning comes from once again thinking of a simple spin as a special case of a frame spin, i.e. as 
$$
[(S^{n+1}, S^{n-1})-S^1\times (B^{n},B^{n-2})]
\bigcup_{\bd} [S^1\times (B_K^{n}, B_K^{n-2})],
$$
where $\bigcup_{\bd}$ indicates gluing along the common boundary in the obvious (untwisted) fashion. Suppose now that we have a $1$-parameter family $f_{\psi}$ of \emph{deformations} of $B^n_K$ rel $\bd B_K^n$ such that $f_0$ is the identity and  $f_{2\pi}(B_K^{n-2})=B_K^{n-2}$.  
The family $f_{\psi}$ should also depend piecewise linearly on the parameter $\psi$. Litherland then describes the deform spin of $K$ as 
$$
[(S^{n+1}, S^{n-1})-S^1\times (B^{n},B^{n-2})]
\bigcup_{\bd} (S^1\times B_K^{n}, \bigcup_{\psi\in S^1}\psi\times f_{\psi}(B_K^{n-2})),
$$
In other words, as we spin, we deform the knotted ball according to  $f_{\psi}$. Note that in this description $ S^1\times B_K^n $ is the ordinary undeformed product, but we equally well could have used the deformation of the pair; Litherland demonstrates the equivalence of the two approaches and uses it to redefine the deform spin in terms of \emph{crossed products} of spaces. However, for our purposes in the following sections, it is perhaps easier to maintain the original viewpoint. In this language, it is easy to observe that simple spinning corresponds to setting $f_{\psi}$ equal to the identity for each $\psi$, while $k$-twist spinning corresponds to setting $f_{\psi}$ equal to rotation of the longitude coordinate of $B_K^n$ by $k\psi$ (technically, to get around the fact that we need to keep the boundary of $B_K^n$ fixed, we rotate a smaller interior ball allowing the region between the two boundaries to become stretched around; however, it is clear that so long as the ball being rotated encompasses the knotted part of $B_K^{n-2}$, this does not affect the final construction of the deform spun knot, and we recover our original description of twist spinning). 

Litherland also shows that, thinking of the collection $f_{\psi}$ as a PL  map $f\co  B_K^n\times [0,2\pi] \to B_K^n\times [0,2\pi]$, the type of the deform spun knot is dependent only upon the \emph{pseudo-isotopy class of $f$ rel $\bd  B_K^n$ as a map of pairs} ($f$ and $g$ are pseudo-isotopic rel $\bd  B_K^n$ as maps  of pairs  if there is a PL homeomorphism $H\co  (B_K^n, B_K^{n-2})\times [0,2\pi] \to (B_K^n, B_K^{n-2})\times [0,2\pi]$ such that $H|_{B_K^n\times 0}$ and $H|_{\bd B_K^n \times [0,2\pi]}$ are the identity maps and $H|_{B_K^n\times 2\pi}=f_{2\pi}g_{2\pi}^{-1}$.)

We can now define roll spinning of a classical knot $K\co  S^1\into S^3$. We present a new geometric description of Litherland's construction that we hope will be valuable to our readers. This interpretation of role spinning is due independently to Dennis Roseman.

 Recall our definition of $(B_K^3, B_K^{1})$ by removing an unknotted ball neighborhood of a point of $K$. Since we are dealing with a classical knot, we can parametrize $S^1$ by angles $\psi$ and consider $(B_{K,\psi}^3, B_{K,\psi}^{1})$ built as the complement of the neighborhood of the point $K(\psi)$ of the knot. We have already noted that for different choices of $\psi$, the pairs $(B_{K,\psi}^3, B_{K,\psi}^{1})$ are all PL homeomorphic. Nevertheless, starting from a fixed base, say $\psi=0$, we can view the collection of homeomorphisms  $f_{\psi}\co (B_{K,0}^3, B_{K,0}^{1})\to (B_{K,\psi}^3, B_{K,\psi}^{1})$ as a one parameter family of deformations and use this to deform spin. This construction is roll spinning. $k$-roll spinning can be created by rolling the basepoint around the knot $k$ times. A more technically precise formulation is given in \cite{Lith} (see also \cite{Ter}).  Note that this construction depends on a choice of framing of $K$ in order to control the roll (in the aeronautical sense) of the ball as it traverses the knot; thus one usually  defines rolling with respect to some fixed standard framing of $K$, usually the one  in which $K$ and a longitude of the boundary torus of the framed neighborhood of $K$ do not link. If we use a different framing, we will twist as we roll; this leads to twist roll-spun knots or, more specifically, $l$-twist $k$-roll spun knots. 

\begin{figure}[!h]
\begin{center}
\scalebox{.5}{\includegraphics{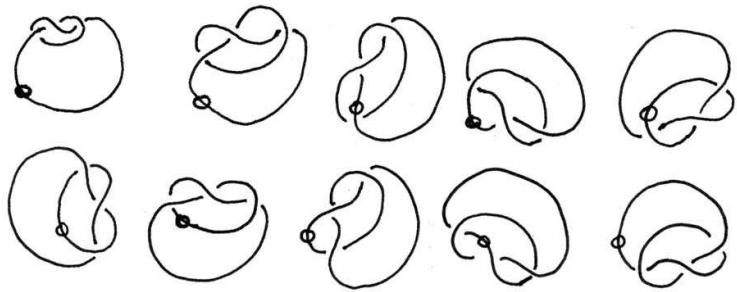}}
\end{center}
\caption{Rolling the trefoil. The circle in these pictures represents $B_{-}$. Rather than moving $B_{-}$, it is more illustrative to hold it fixed and roll the knot around it!  In fact, this roll of the trefoil is twisted with respect to the unlinked longitude. Furthermore, due to the symmetries of a torus knot, the rolling part of the deformation has no effect, only the twisting. Can you see why? Nonetheless, this diagram illustrates the idea of the procedure.}
\end{figure}

Litherland treats another example of deform spinning that applies only to knots which possess symmetries, i.e. periodic homeomorphisms $S^n\to S^n$ that take the knot to itself. This construction is called \emph{symmetry spinning}. 
However, the construction is slightly technical, involving certain branched covers of the sphere, and so we omit a detailed description. The interested reader is referred to \cite{Lith} or \cite{Kan}, in which Taizo Kanenobu utilizes symmetry spinning to obtain some counterexamples to the $4$-dimensional Smith conjecture.

Unfortunately, roll-spinning cannot be generalized to roll spinning of classical knots about higher dimensional spheres. As for twist-spinning, this is because $\pi_n(S^1)=0$ for  $n>1$, so any attempt to parametrize rolling by a sphere of dimension $>1$ will re-create superspinning.  Similarly, since the $n$th homotopy group $\pi_n(S^1\times S^1)$ is trivial for $n>1$,   super twist-roll spinning is also nothing new. 
However, Litherland does discuss generalizations such as roll spinning  higher-dimensional knots about  $S^1$. This utilizes an isotopy of $S^{n-1}$ just as classical roll-spinning utilizes the isotopy that rotates $S^1$. Rather than visit the technical details here, we will jump right to a more general construction.

\subsection{Frame deform spinning}\label{S: frame def}

Putting together frame spinning and deform spinning, we can introduce a new knot construction, \emph{frame deform spinning} (this construction is implicit in a remark in \cite[\S 3]{Ro89}). By now the method may  be obvious: we begin with an $n$-knot $K$ and an $m$-manifold $M^m$ embedded with framing in $S^{n+m-2}$, which itself sits unknotted and with the standard framing in $S^{n+m}$. We also posit a map $f$ from $M^m$ into the space of PL homeomorphisms of $B_K^n$ rel $\bd B_K^n$, taking $x\in M^m$ to $f_x$.  Then we can define the frame deform spin of $K$ as 
$$
[(S^{n+m}, S^{n+m-2})-M^m\times (B^{n},B^{n-2})]
\bigcup_{M^m\times S^{n-3}} (M^m\times B_ K^{n}, \cup_{x\in M^m} f_x(B_K^{n-2})).
$$
If $K$ is a classical knot  and there is a non-trivial map $g\co  M^m\to S^1$, we can compose $g$ with the $1$-parameter family of deformations used to define roll-spinning to create \emph{frame roll spinning}. We could also use a map $M\to S^1\times S^1$ to define \emph{frame twist-roll spinning}.

An example of frame deform spinning due to Roseman provides a higher-dimensional version of rolling. In classical rolling, we created a $1$-parameter deformation by moving the point  at which we cut out $(B^3_-, B^{1}_-)$ around the knotted circle $K$. What if we similarly moved the excision point around a higher-dimensional knot in order to roll spin a knotted $S^p$ about $S^p$? The difficulty is that as we move the point around, we really need a frame on it to describe precisely how the ball $B^n_-$ is situated in $S^n$ so that we can say exactly how the complement is deformed under the motion. For non-deform spinning, the frame is irrelevant since all are equivalent by a homeomorphism and so yield the same spun knots. For roll spinning, however, the choice is relevant as our deformation depends directly on the global movement of the frame. In classic rolling, this frame was completely determined by our previous considerations -  by the amount of twisting orthogonal to the knot and by the direction of the rolling (positive or negative). For higher dimensions, our previous considerations eliminate twisting and allow us to find a trivial normal framing, orthogonal to the knot, but the framing parallel to the knot will exist at all points only if the sphere is parallelizable. Thus this type of rolling is only possibly if $p$ equals $1$, $3$, or $7$, yielding knots of dimension $4$, $8$, or $16$. Very little is known as yet about this construction and the knots it generates.

\section{Other constructions}

We close by briefly mentioning two other known constructions of knots related to spinning. 

The first, due to John Klein and Alexander Suciu \cite{KS91}, is called \emph{diff-spinning}. It is a modified version of frame spinning in the smooth category in which the manifold $M^m$ is altered by a diffeomorphism in the process of spinning. Klein and Suciu used diff-spinning to demonstrate the existence of  inequivalent fibered knots whose homotopy Seifert pairings are isometric. It had been shown by Michael Farber in \cite{Fa80} that  the homotopy Seifert pairing uniquely determines a fibered knot, provided the knot is sufficiently simple (i.e. its complement has the homotopy type of $S^1$ in sufficiently many dimensions). The construction of Klein and Suciu showed that Farber's result does not extend for all knots. 

To define diff-spinning, note that the complement of the  frame spin about $M^m$ is diffeomorphic
to the union of $B^{n+m-1}\times S^1$ with $M^m \times X$, where $X$ is the complement of $K$ (see \cite{KS91} for a diagram illustrating this). 
Suppose now that we are given a self-diffeomorphism $\Phi$ of $M^m$ that extends to a diffeomorphism $\bar \Phi$ of $B^{n+m-1}\supset S^{n+m-2}\supset M^m$.    Then, roughly speaking, the diff-spin is formed by removing this complement and replacing it with the twisted product $(B^{n+m-1}\times_{\bar \Phi} S^1)\bigcup(M^m\times_{\Phi}X)$. If $\Phi$ satisfies a certain algebraic condition (see \cite[\S 5]{KS91}), this space will also be the complement of a knot, the diff-spun knot. 

Spinning in the smooth category is also considered in a low-dimensional context by Ronald Fintushel and Ronald J. Stern in \cite{FS97}. Their construction, called \emph{rim surgery}, begins with a surface $\Sigma$ embedded in a smooth $4$-manifold $X$. The procedure then is to spin a classical knot $K\subset S^3$ around a certain curve embedded in $\Sigma$. The technical details are essentially those we have seen before: a trivial neighborhood pair of the curve is removed, and it is replaced with a neighborhood pair whose fibers are $(B^3_K, B^1_K)$. The result is a new surface embedding $(X,\Sigma_K)$. Under certain assumptions, $(X,\Sigma_K)$ will be homeomorphic to $(X,\Sigma)$, but if the Alexander polynomial of $K$ is nontrivial, $(X,\Sigma_K)$ and $(X,\Sigma)$ will not be diffeomorphic! More generally, $(X,\Sigma_{K_1})$ and $(X,\Sigma_{K_2})$ will be diffeomorphic only if the two knots $K_1$ and $K_2$ have the same Alexander polynomial. This is proven as an application of Seiberg-Witten theory. Generalizations of this construction and various applications have been considered by a variety of authors, including Fintushel and Stern  \cite{FS99}, Stefano Vidussi \cite{V1,V2}, Sergey Finashin \cite{Fi02}, and Hee Jung Kim \cite{Kim}, who considers a twist-spinning analogue. 

One last type of spinning, also introduced by Roseman in \cite{Ro89}, is what he calls ``spinning a knot about a projection''. We shall refer to this as \emph{projection spinning}. This clever construction involves many technical details, but, very roughly, the idea is to spin about an immersed manifold $M$, rather than an embedded one as we did for frame spinning. Let $\Sigma$ denote the singular set of the embedding, i.e. the image of the points of $M$ for which the immersion is not $1$-$1$. Let $N(\Sigma)$ be a neighborhood of $\Sigma$ in $S^n$. Outside of $N(\Sigma)$, the construction is the same as for frame spinning, i.e. each unknotted fiber pair in the tubular neighborhood of $M-N(\Sigma)$ is replaced with $(B^n_K, B^{n-2}_K)$. However, $N(\Sigma)$ itself is broken up into sets  that are homeomorphic to  balls $B^n$ and such that the intersection of $B^n$ with the immersed $M$ is a collection of hyperplanes. These balls with hyperplanes are then used as the data to create \emph{multi-knots}, in which some knot $K$ is blended together with itself in multiple directions. These ball neighborhoods are then replaced with the multi-knots. If $M$ is embedded, we recover frame-spinning.  We refer the reader to \cite{Ro89} both for the technical definitions of projection spinning and for some nice graphical illustrations of the process. 

In Remark 7 of \cite{Ro89}, Roseman notes that it is further possible to deform projection spin, perhaps the ultimate in knot spinning constructions!

\bibliographystyle{amsplain}
\bibliography{bib}

\end{document}